\documentclass[psamsfonts]{amsart}

\usepackage{amssymb,amsfonts,amsthm}
\usepackage[all,arc]{xy}
\usepackage{enumerate}
\usepackage{mathrsfs}
\usepackage{hyperref}
\usepackage{dsfont}
\hypersetup{
pdftitle={Paper Xp},
pdfsubject={Mathematics, Harmonic Analysis, Functional Analysis},
pdfauthor={},
pdfkeywords={}
}
\usepackage{bm}
\usepackage{pgfplots}

\allowdisplaybreaks

\usepackage{enumitem, hyperref}
\makeatletter
\def\namedlabel#1#2{\begingroup
    #2%
    \def\@currentlabel{#2}%
    \phantomsection\label{#1}\endgroup
}

\newtheorem{theorem}{Theorem}[section]
\newtheorem{corollary}[theorem]{Corollary}
\newtheorem{problem}[theorem]{Problem}
\newtheorem{proposition}[theorem]{Proposition}
\newtheorem{lemma}[theorem]{Lemma}

\newtheorem{remark}[theorem]{Remark}
\newtheorem*{remark*}{Remark}

\numberwithin{equation}{section}
\numberwithin{theorem}{section}

\theoremstyle{theorem}
\newtheorem{ltheorem}{Theorem}

\newcommand{\C}{\mathbb{C}}
\newcommand{\Z}{\mathbb{Z}}
\newcommand{\N}{\mathbb{N}}

\newcommand{\suma}[1]{\sum\limits_{#1}}
\newcommand{\Suma}[2]{\sum\limits_{#1}^{#2}}

\newcommand*{\bigchi}{\mbox{\large$\chi$}}

\newcommand{\dem}{\noindent {\bf Proof. }}
\newcommand{\demA}{\noindent {\bf Proof of Theorem A. }}
\newcommand{\fin}{\hspace*{\fill} $\square$ \vskip0.2cm}

\newcommand{\cexp}{\mathsf{E}}
\newcommand{\mult}{\mathrm{M}}
\newcommand{\symbolmult}{\mathrm{m}}
\newcommand{\Sb}{\mathsf{S}}

\newcommand{\Xp}{\mathrm{X}_p}
\newcommand{\Hgroup}{\mathrm{H}}
\newcommand{\HGroup}{\Gamma}
\newcommand{\Ggroup}{\mathrm{G}}
\newcommand{\sgn}{\mathrm{sgn}}

\title[Metric $\mathrm{X}_p$ on groups]{Trigonometric chaos and $\mathrm{X}_p$ inequalities II: \\ $\mathrm{X}_p$ inequalities with sharp scaling parameter}

\begin{document}

\null

\vskip-15pt

\null

\begin{center}
{\huge Trigonometric chaos and $\mathrm{X}_p$ inequalities II} \\ {\Large ---$\mathrm{X}_p$ inequalities in group von Neumann algebras---}

\vskip15pt

{\sc {Antonio Ismael Cano-M\'armol, \\ Jos\'e M. Conde-Alonso and Javier Parcet}}
\end{center}

\title[]{}

\maketitle

\null

\vskip-60pt

\null

\begin{center}
{\large {\bf Abstract}}
\end{center}

\vskip-30pt

\null

\begin{abstract}
In the line of previous work by Naor, we establish new forms of metric $\Xp$ inequalities in group algebras under very general assumptions. Our results' applicability goes beyond the previously known setting in two directions. In first place, we find continuous forms of the $\mathrm{X}_p$ inequality in the $n$-dimensional torus. Second, we consider transferred forms of the sharp scalar valued metric $\mathrm{X}_p$ inequality in the von Neumann algebra $\mathcal{L}(\Ggroup)$ of a discrete group $\Ggroup$. As a byproduct of our results, some metric consequences and their relation with bi-Lipschitz nonembeddability of Banach spaces are explored in the context of noncommutative $L_p$ spaces.   
\end{abstract}

\addtolength{\parskip}{+1ex}

\vskip15pt

\section*{\bf \large Introduction}

In functional analysis, it is a fundamental question to know when a given space is isomorphic to a linear subspace of another. The corresponding area, Banach space embedding theory, dates back to the first half of the twentieth century. For the important example of $L_p$ and $L_q$ the landscape is well known. $L_2$ is isomorphic to a subspace of $L_p$ for all $p$ in the Banach range, but there is no linear embedding from $L_q$ to $L_p$ if either $q<\min\{2,p\}$ or $q>\max\{2,p\}$. Banach conjectured a positive answer for $\min\{2,p\} < q < \max\{2,p\}$. Kadec proved it for $p<q<2$ in \cite{Ka}, while Paley disproved it for $2<q<p$ in \cite{P36}. 

Sometimes, properties of Banach spaces that are local (that is, ones that can be formulated in terms of vectors in finite dimensional subspaces) do not rely on the whole strength of the linear structure. Instead, they may only need the fact that the norm defines a metric. Determining for which properties is this the case is known as the Ribe program, initiated after \cite{R76} and explicitly formulated in \cite{B86}. Replacing the linear maps that yield inclusions between Banach spaces by bi-Lipschitz maps becomes therefore natural. Regarding $L_p$ spaces, the qualitative results mentioned above do not change, and notably one gets again that $L_q$ cannot be embedded into $L_p$ in the range $2<q<p$. This happens in spite of the fact that the classical argument that shows that there is no bi-Lipschitz embedding from $L_q$ into $L_p$ uses explicitly the linear structure of the spaces. This is so because it goes through differentiability |of Lipschitz maps| to reduce the metric statement to a linear one. However, the result can be proved by other means and indeed a different quantitative proof can be found in the work of Naor and Schechtman \cite{NS16}. One can consult the Introduction in that work for a more detailed context and references regarding the history of the problem and its connections with other areas. 

A key innovation in \cite{NS16} is the introduction of a family of inequalities, termed $\Xp$ inequalities. Let us briefly recall the results in \cite{NS16}, since they are going to be the starting point of our investigation in the present manuscript. Denote by $[n]$ the initial segment $[n]:= \{1,2, \ldots , n\}$. We view the set $\mathbb{Z}_{8m}^n$ as a probability space equipped with the normalized counting measure. In the same vein we identify $\mathbb{Z}_2^n$ with $\{-1,1\}^n$ using multiplicative notation. A metric space $(\mathbb{X},d)$ is said to be a metric $\mathrm{X}_p$ space if for each $n\in\mathbb{N}$ and $k\in[n]$, there exists $m\in\mathbb{N}$ such that every mapping $f:\mathbb{Z}_{8m}^n \to \mathbb{X}$ satisfies the following estimate for $p \ge 2$
\begin{align}
\frac{1}{{n \choose k}} \suma{\substack{\Sb \subset [n] \\ |\Sb|=k}} \cexp & \big[ d \big( \mult_{4m\varepsilon_{\Sb}}f(x),f(x) \big)^p \big] \nonumber \\ & \lesssim_p m^p \left( \frac{k}{n} \sum_{j=1}^n \cexp \big[ d \big( \mult_{e_j}f(x),f(x) \big)^p\big] +\Big(\frac{k}{n}\Big)^{\frac{p}{2}}\cexp \big[ d \big( \mult_{\varepsilon}f(x),f(x) \big)^p \big] \right). \label{eq:metricXp} \tag{$\mathrm{MX}_p$}
\end{align}
Above, the dummy variable $(\varepsilon,x)$ belongs to $\mathbb{Z}_2^n \times \mathbb{Z}_{8m}^n$, $\varepsilon_{\Sb}$ denotes the $\Sb$-truncated vector $\sum_{j \in \Sb} \varepsilon_j e_j$, the expectation $\cexp$ is taken in the product probability space, and the operators $\mult_{v}$ are just translations given by
$$
\mult_v f(x) = f(x+v).
$$
They are Fourier multipliers, which will become relevant later. The main result in \cite{NS16} then asserts that $L_p$ is an $\Xp$ space, while $L_q$ is not if $2<q\neq p$. This and the fact that being an $\mathrm{X}_p$ metric space is invariant under the action of bi-Lipschitz maps provides the quantitative nonembeddability argument alluded above. Note that the result is not true independently of the value of $m$. It is necessary that it scales at least so that $m \gtrsim \sqrt{n/k}$, even though the bound was not found yet to be sharp in \cite{NS16}. The above discussion justifies why $\Xp$ inequalities are a tool to study Lipschitz nonembeddability results between Banach spaces.

Little after it, Naor found in \cite{N16} a different proof of the metric $\Xp$ inequality which was quantitatively optimal. Namely, it turns out that the necessary bound $m \gtrsim \sqrt{n/k}$ is, indeed, sufficient. This entails a number of significant consequences in metric geometry, including distortion estimates of discrete $\ell_q$-grids into $L_p$ and sharp snowflake $L_p$-estimates. A key step towards the main result in \cite{N16} is the inequality below for balanced averages of Fourier truncations of the Fourier-Walsh expansions. Those are expansions of the form
$$
f= \sum_{\mathsf{A} \subset [n]} \widehat{f}(\mathsf{A}) W_\mathsf{A}.
$$ 
Given $2\leq p <\infty$, $k \in [n]$, and $f \in L_p(\mathbb{Z}_2^n)$ mean zero, it holds 
\begin{equation}\label{eq:Naor} \tag{$\mathrm{N}_p$}
\frac{1}{{n \choose k}} \sum_{\substack{\Sb \subset [n] \\ |\Sb|=k}} \Big\| \sum_{\mathsf{A} \subset \mathsf{S}} \widehat{f}(\mathsf{A}) W_\mathsf{A} \Big\|_p^p \lesssim_p \frac{k}{n} \sum_{j=1}^n \left\| \partial_j f\right\|_p^p + \Big( \frac{k}{n} \Big)^{\frac{p}{2}} \|f\|_p^p.
\end{equation}
Above, the symbol $\partial_j$ denotes the $j$-th directional discrete derivative
$$
\partial_j f(\varepsilon) := f(\varepsilon) - f(\varepsilon_1, \ldots, \varepsilon_{j-1}, - \varepsilon_j, \varepsilon_{j+1}, \ldots, \varepsilon_{n}).
$$
Once \eqref{eq:Naor} has been established, the metric $\Xp$ inequality for all classical $L_p$ spaces follows from it plus the natural injection of the multiplicative group $\mathbb{Z}_2^n$ into $\mathbb{Z}_{2m}^n$ needed to define the multipliers $\mult_\varepsilon(f) (x) = f(x+\varepsilon)$ for $f:\mathbb{Z}_{2m}^n \to \mathbb{C}$. This is the starting point of our work. We seek new metric $\Xp$ inequalities that can be proven using appropriate versions of \eqref{eq:Naor}. In \cite{CCP22}, a wide range of inequalities of this kind is established for the (noncommutative) $L_p$ spaces associated with the group von Neumann algebra $\mathcal{L}(\Ggroup)$ of a discrete group $\Ggroup$. Elements in $L_p(\mathcal{L}(\Ggroup))$ are not functions for nonabelian $\Ggroup$ and making sense of pointwise evaluation is not immediate, though necessary. Seeking the greatest possible generality, we shall consider below different choices of abelian groups $\Hgroup$ that will play the role of $\mathbb{Z}_2^n$ in  \cite{N16}. Therefore, in the evaluation $\mult_\varepsilon f(x):= f(x+\varepsilon)$, the variable $\varepsilon \in \Hgroup$ but that may not be the case for $x$, which imposes a compatibility condition between the two variables. This is a crucial difficulty and the main contribution in this paper is to clarify when a chaos-type inequality like \eqref{eq:Naor} |as those obtained in \cite{CCP22}| can be used to deduce a new metric $\Xp$-type inequality. 

Let us now explain our abstract result, we refer to \cite{CCP22} for notation regarding analysis on group von Neumann algebras. Let $\Hgroup$ be a discrete abelian group which can be written as a direct product:
$$
\Hgroup = \Hgroup_1 \times \Hgroup_2 \times \cdots \times \Hgroup_n, 
$$ 
equipped with its counting measure. Given $\Sb \subset [n]$, we define the truncation of $\Hgroup$ on $\Sb$ as the subgroup $\Hgroup_\Sb = \times_{j \in \Sb} \Hgroup_j \hookrightarrow \Hgroup$. Truncations are transferred to the dual group $\widehat{\Hgroup}$, which is a probability space with Haar measure $dx$. 
Assume also that $\Ggroup$ is an arbitrary discrete group, and let $p\geq 2$. We say that $(\Hgroup,\Ggroup)$ is an $\Xp$\emph{-representable pair} if the following conditions hold:
\begin{itemize}
\item[\namedlabel{itm:2}{(RP-1)}] \textbf{Balanced truncations.} We have 
$$\hskip20pt \frac{1}{{n \choose k}} \suma{\substack{\Sb \subset [n] \\ |\Sb|=k}} \big\|(\cexp_{[n]\setminus \Sb} \otimes \mathrm{id}) h\big\|^p_p \lesssim_p \frac{k}{n} \Suma{j=1}{n} \big\|(\partial_j \otimes \mathrm{id}) h \big\|^p_p + \Big(\frac{k}{n}\Big)^{\frac{p}{2}} \|h\|^p_p,$$
for $k \in [n]$, any mean zero $h \in L_p(\widehat{\mathrm{H}};L_p(\mathcal{L}(\mathrm{G})))$, and where 
$$
\hskip20pt \cexp_{[n]\setminus \Sb}f(x) = \int_{\widehat{\Hgroup}} f(x_{\Sb} + z_{[n]\setminus \Sb}) \ dz
$$
is the conditional expectation onto functions which only depend on variables in $\widehat{\Hgroup}_\Sb$. The $\partial_j$ are linear maps which play the role of directional derivatives.

\vskip8pt

\item[\namedlabel{itm:12}{(RP-2)}] \textbf{Compatibility conditions.} There exists an abelian group $$\hskip20pt \HGroup =  \HGroup_1 \times \HGroup_2 \times \cdots \times \HGroup_n$$ and a map $\eta: \widehat{\Hgroup} \rightarrow \widehat{\HGroup}$ such that the following compatibility conditions hold:
\vskip2pt
\begin{itemize}
\item[i)] \textbf{Symmetric inclusion.} Given $\Sb \subset [n]$, the variables $\eta_\Sb(y):=(\eta(y))_\Sb$ and $-\eta_\Sb(y)$ (using here additive notation) have the same distribution.
\vskip2pt
\item[ii)] \textbf{Uniformly bounded translations.} There exists a group of Fourier multipliers $\{\mult_\gamma \}_{\gamma \in  \widehat{\HGroup}}$, completely and uniformly bounded on $L_p(\mathcal{L}(\Ggroup))$.
\end{itemize}
\end{itemize}

\begin{ltheorem}\label{th:abstractThm}
Let $p \geq 2$ and $(\Hgroup,\Ggroup)$ be an $\Xp$-representable pair. Let $k \in [n]$ and suppose that $m \in \N$ satisfies $m \geq \sqrt{n/k}$. Then, for every $f \in L_p(\mathcal{L}(\Ggroup))$ there holds
	\begin{align*}
		\frac{1}{{n \choose k}} &\sum_{\substack{\Sb \subset [n] \\ |\Sb|=k}} \int_{\widehat{\Hgroup}} \big\|\mult_{4m\eta_\Sb(y)}f - f\big\|^p_{p} \ dy \\
		&\lesssim_p m^p \Bigg(\frac{k}{n} \Suma{j=1}{n} \int_{\widehat{\Hgroup}} \big\|\partial_j \mult_{2 \eta(y)}f  
		\big\|_p^p \ dy + \Big(\frac{k}{n} \Big)^{\frac{p}{2}} \int_{\widehat{\Hgroup}} \big\|\mult_{\eta(y)}f - f \big\|_p^p \ dy \Bigg).
\end{align*}
\end{ltheorem}

Metric $\mathrm{X}_p$ inequalities from \cite{N16,NS16} apply for functions $f:\mathbb{Z}_{8m}^n \to \mathbb{X}$. Theorem \ref{th:abstractThm} is an inequality for $f \in L_p(\mathcal{L}(\Ggroup))$, so that it should be regarded as a generalization where we replace $\mathbb{Z}_{8m}^n$ by the dual group of $\Ggroup$. Note that its dual is a quantum group rather than a classical group when $\Ggroup$ is nonabelian. Additionally, the role of $\Hgroup$ is merely to provide an index set over which translations are performed, which generalizes the former role of $\mathbb{Z}_2^n$. In particular, Theorem \ref{th:abstractThm} generalizes  \eqref{eq:metricXp} with $(\Hgroup,\Ggroup) = (\mathbb{Z}_{2}^n,\mathbb{Z}_{8m}^n)$, auxiliary group $\Gamma = \mathbb{Z}_{8m}^n$ and $\eta$ being the natural inclusion defined coordinatewise by $1 \mapsto 1$ and $(-1) \mapsto 8m-1$. Note that 
\begin{equation} \tag{$\partial_j$} \label{EqDrivative}
\int_{\mathbb{Z}_2^n} \big\|\partial_j \mult_{2\eta(y)} f
\big\|_{p}^p \ dy \lesssim_p  \|\mult_{e_j}f - f \|_{p}^p
\end{equation} 
holds when $\partial_j$ is the aforementioned discrete derivative in the $j$-th direction. The role of $\HGroup$ was hidden in \cite{N16,NS16}. This is natural since it is possible to formulate Theorem \ref{th:abstractThm} without appealing to the existence of $\eta$ and $\HGroup$. Indeed, one can find axioms on a family of multipliers $\{\mult_y\}_{y \in \Hgroup}$ and their powers |so that they act like translations| which are formally weaker than \ref{itm:2} and \ref{itm:12} and suffice to get the conclusion of Theorem \ref{th:abstractThm}. We have chosen to formulate our result in this manner because in all of our examples we may identify the auxiliary objects $\eta$ and $\HGroup$, which is an efficient way of checking that we can apply Theorem \ref{th:abstractThm}. However the existence of weaker conditions can be useful for future constructions and we will detail them in Section \ref{sec1}. In the general case, the improved statement that one gets when an inequality like \eqref{EqDrivative} holds is a true metric inequality, while the conclusion of Theorem \ref{th:abstractThm} is not in general. In our other examples, we shall choose derivatives whose relation with the multipliers $\mult_v$ yields the stronger, metric statement. 

The proof of Theorem \ref{th:abstractThm} follows the strategy in \cite{N16}. Our main contribution is the identification of the right conditions under which said strategy can be extended to numerous contexts, including noncommutative ones. We illustrate Theorem \ref{th:abstractThm} with several concrete scenarios, some of which we describe next. In each case, \ref{itm:2} is checked using the results in \cite{CCP22} and we focus in getting fully metric statements.
\begin{itemize}
\item[i)] \textbf{Continuous metric $\Xp$ inequalities.} Taking $\Hgroup=\Ggroup=\mathbb{Z}^n$, we obtain two continuous statements in the $n$-dimensional torus using two different choices for the inclusion $\eta$. One of them can be obtained from  \cite{N16} by an approximation procedure, while the other seems to require Theorem \ref{th:abstractThm}.

\vskip5pt

\item[ii)] \textbf{Cyclic extensions of metric $\Xp$ inequalities.} Let $\eta_{\ell}$ be the inclusion that maps $\mathbb{Z}_{2\ell}^n$ into $\{k: \; k \in [\ell] \mbox{ or } 8m\ell -k \in [\ell] \} \subset \mathbb{Z}_{8m\ell}^n$. We prove a metric $\mathrm{X}_p$ inequality which reduces to \eqref{eq:metricXp} for $\ell=1$. It is sharp in terms of the scaling parameter $m$ and entails the same metric consequences as \cite{N16}. 

\vskip5pt

\item[ii)] \textbf{Transferred $\Xp$ inequalities.} When $\Ggroup$ is nonabelian, we first observe that our results in \cite{CCP22} hold with values on noncommutative $L_p$ spaces over QWEP algebras, so that \ref{itm:2} is satisfied. Next, using an appropriate corepresentation for $(\Hgroup,\Ggroup)$ we obtain $\mathrm{X}_p$ inequalities on the free group $\mathbb{F}_n$. 
\end{itemize}

The above scenarios show the applicability of our main result. In metric terms we observe that noncommutative $L_p$ spaces are metric $\Xp$ spaces. This implies that the distortion estimates for $\ell_q$-grids and $L_p$-snowflakes in \cite{N16} still hold when considering embeddings into noncommutative $L_p$. Beyond this, we have not found significantly new applications of our main results in the context of metric geometry so far. It would be very interesting to generalize Theorem \ref{th:abstractThm} so as to include $\Xp$-representable pairs $(\Hgroup,\Ggroup)$ admitting \emph{nonabelian translations} in the index set $\Hgroup$.



\section{\bf \large Metric $\Xp$ inequalities} \label{sec1}

We work in the language of noncommutative $L_p$ spaces over group von Neumann algebras, as explained in \cite{CCP22}. Given a discrete group $\Ggroup$, its von Neumann algebra is denoted by $\mathcal{L}(\Ggroup)$ |which coincides with $L_\infty(\widehat{\Ggroup})$ when $\Ggroup$ is abelian|, while $L_p(\mathcal{L}(\Ggroup))$ is the associated noncommutative $L_p$ space, $2\leq p < \infty$. $L_p^\circ(\mathcal{L}(\Ggroup))$ will be the family of zero-trace elements in $L_p(\mathcal{L}(\Ggroup))$. Throughout this section assume that $n$ is fixed. We emphasize again that we follow the strategy of proof in \cite{N16} for an $\Xp$-representable pair $(\Hgroup,\Ggroup)$, as defined in the Introduction. Given $\Sb \subset [n]$ consider the auxiliary operator $T_\Sb$ on $\mathcal{L}(\Ggroup)$ given by 
\begin{align*}
    T_\Sb f = \int_{\widehat{\Hgroup}} \mult_{2 \eta_\Sb(y)} f \ dy, \quad \mbox{for} \quad f \in L_p(\mathcal{L}(\Ggroup)).
\end{align*}

\demA
For each $\Sb \subset [n]$ with $|\Sb|=k$ and $y \in \Hgroup$
\begin{eqnarray*}
\big\| \mult_{4m\eta_\Sb(y)}f - f \big\|_p^p & \lesssim_p & \big\| T_{[n]\setminus \Sb}f - f \big\|_p^p \\ & + & \big\| \mult_{4m\eta_\Sb(y)}T_{[n]\setminus \Sb}f - T_{[n]\setminus \Sb}f \big\|_p^p \\ & + &  \big\| \mult_{4m \eta_\Sb(y)}f - \mult_{4m \eta_\Sb(y)}T_{[n]\setminus \Sb}f \big\|_p^p \ =: \ \mathrm{A}+\mathrm{B}+\mathrm{C}.
\end{eqnarray*}
First, we claim that 
\begin{equation} \label{lema8}
\big\| T_\Sb f - f \big\|_p^p \lesssim_p \int_{\widehat{\Hgroup}} \big\| \mult_{\eta(y)}f - f \big\|_p^p \, dy.
\end{equation}
Indeed, by convexity of the $L_p$ norm
\begin{eqnarray*}
        \big\|T_\Sb f - f \big\|_p^p & = & \Big\| \int_{\widehat{\Hgroup}} \mult_{2 \eta_\Sb(y)}f \, dy - f \Big\|^p_p \ \lesssim_p \ \int_{\widehat{\Hgroup}} \big\| \mult_{2 \eta_\Sb(y)}f - f \big\|_p^p \, dy \\
        & \leq & 2^{p-1} \left( \int_{\widehat{\Hgroup}} \big\| \mult_{2 \eta_\Sb(y)}f - \mult_{\eta(y)}f \big\|_p^p \, dy+  \int_{\widehat{\Hgroup}} \big\| \mult_{\eta(y)}f - f \big\|_p^p \, dy \right).
    \end{eqnarray*}
Now, $\eta(y)= \eta_S(y) + \eta_{[n]\setminus S}(y)$ and $\eta_S(y) - \eta_{[n]\setminus S}(y)$ are identically distributed by \ref{itm:12}. Therefore, by the properties of the family $\{\mult_{\varepsilon}\}_\varepsilon$ |also in \ref{itm:12}| the first term above satisfies
\begin{align*}
        \int_{\widehat{\Hgroup}} \|\mult_{2 \eta_\Sb(y)}f - \mult_{\eta(y)}f \|_p^p \ dy & = \int_{\widehat{\Hgroup}} \|\mult_{\eta(y)} \mult_{\eta_\Sb(y)-\eta_{[n]\setminus \Sb}(y)}f - \mult_{\eta(y)}f\|_p^p \ dy \\
        &\lesssim \int_{\widehat{\Hgroup}} \|\mult_{\eta_\Sb(y)-\eta_{[n]\setminus \Sb}(y)}f - f\|_p^p \ dy \\
        &= \int_{\widehat{\Hgroup}} \|\mult_{\eta(y)}f - f\|_p^p \ dy,
\end{align*}
which implies the claim. This proves that
$$
\mathrm{A}+\mathrm{C} \, \lesssim \, \mathrm{A} \, \lesssim \, \int_{\widehat{\Hgroup}} \big\| \mult_{\eta(y)}f - f \big\|^p_p \, dy.
$$
For the term $\mathrm{B}$, our claim is the following
\begin{eqnarray} \label{lema9}
\lefteqn{\hskip-35pt \frac{m^{-p}}{{n \choose k}} \suma{\substack{\Sb \subset [n] \\ |\Sb|=k}} \int_{\widehat{\Hgroup}} \big\| \mult_{4m \eta_\Sb(y)} T_{[n]\setminus \Sb}f - T_{[n]\setminus \Sb}f \big\|_p^p \, dy} \\ \nonumber \hskip40pt \!\!\!\! & \lesssim_p & \!\!\!\! \frac{k}{n} \Suma{j=1}{n} \int_{\widehat{\Hgroup}} \big\| \partial_j \mult_{2 \eta(y)}f  \big\|^p_p \, dy + \Big(\frac{k}{n}\Big)^{\frac{p}{2}} \int_{\widehat{\Hgroup}} \big\| \mult_{\eta(y)}f - f \big\|^p_{p} \, dy.
\end{eqnarray}
To prove \eqref{lema9}, we start once more with \ref{itm:12} and the triangle inequality
\begin{eqnarray} \label{eq:lemma2Eq1}
\lefteqn{\Big( \int_{\widehat{\Hgroup}}  \big\| \mult_{4m \eta_\Sb(y)} T_{[n]\setminus \Sb}f - T_{[n]\setminus \Sb}f \big\|_p^p \, dy  \Big)^{\frac{1}{p}}} \\ \nonumber & \leq & \Suma{j=1}{m} \Big( \int_{\widehat{\Hgroup}} \big\| \mult_{4j\eta_\Sb(y)} T_{[n]\setminus \Sb}f - \mult_{4(j-1)\eta_\Sb(y)} T_{[n]\setminus \Sb}f \big\|^p_p \, dy \Big)^{\frac{1}{p}} \nonumber \\ \nonumber & = & \Suma{j=1}{m} \Big( \int_{\widehat{\Hgroup}} \big\| \mult_{(4j-2) \eta_\Sb(y)} [\mult_{2 \eta_\Sb(y)} T_{[n]\setminus \Sb}f - \mult_{-2 \eta_\Sb(y)}T_{[n]\setminus \Sb}f] \big\|_p^p \, dy \ \Big)^{\frac{1}{p}} \nonumber \\ \nonumber & \lesssim & m \ \Big( \int_{\widehat{\Hgroup}} \big\| \underbrace{\mult_{2 \eta_\Sb(y)}T_{[n]\setminus \Sb}f - \mult_{-2\eta_\Sb(y)} T_{[n]\setminus \Sb}f}_{F_\Sb(y)} \big\|_p^p \, dy \ \Big)^{\frac{1}{p}}.
\end{eqnarray}
We claim that $F_\Sb(y) = \cexp_{[n]\setminus \Sb}h(y)$ with $$h: \widehat{\Hgroup} \rightarrow L_p(\mathcal{L}(\Ggroup)),$$ $$h(y) = \mult_{2\eta(y)}f - \mult_{-2\eta(y)}f.$$
Indeed, recall that 
$$
\cexp_{[n]\setminus \Sb}g(x) = \int_{\widehat{\Hgroup}} g(x_{\Sb} + z_{[n]\setminus \Sb}) \, dz.
$$
Also, denote the symbol of $\mult_{\varepsilon}$ by $\symbolmult_{\varepsilon}$, so that
$$
\mult_{\varepsilon} f = \sum_{w \in \Ggroup} \symbolmult_{\varepsilon}(w) \widehat{f}(w) \lambda(w),
$$
where $\lambda(g)$ denotes the left regular representation of $\Ggroup$. Then, given that $\cexp_{[n]\setminus \Sb}$ is linear, and using again the equidistribution property in \ref{itm:12}, we deduce that
	\begin{eqnarray*}
		\cexp_{[n]\setminus \Sb}h(y) & = & \suma{w \in \Ggroup} \widehat{f}(w) \lambda(w) \ [\cexp_{[n]\setminus \Sb}\symbolmult_{2 \eta(y)}(w) - \cexp_{[n]\setminus \Sb}\symbolmult_{-2\eta(y)}(w)]
		\\ & = & \suma{w \in \Ggroup} \widehat{f}(w) \lambda(w) \ [\symbolmult_{2 \eta_\Sb(y)}(w) - \symbolmult_{-2\eta_\Sb(y)}(w)] \ \int_{\widehat{\Hgroup}} \symbolmult_{2 \eta_{[n]\setminus \Sb}(z)}(w) \ dz \\ [5pt]
		& = & \mult_{2\eta_\Sb(y)}T_{[n]\setminus \Sb}f -\mult_{-2\eta_\Sb(y)}T_{[n]\setminus \Sb}f \ = \ F_\Sb(y).
	\end{eqnarray*}	
Moreover, we have 
\begin{align*}
		\int_{\widehat{\Hgroup}} h(y) \ dy & = \int_{\widehat{\Hgroup}} \cexp_{[n]\setminus \Sb}h(y) \ dy \\
		& = \int_{\widehat{\Hgroup}} \big[ \mult_{2 \eta_\Sb(y)}T_{[n]\setminus \Sb}f - \mult_{-2 \eta_\Sb(y)}T_{[n]\setminus \Sb}f \big] \ dy \\
		& = \int_{\widehat{\Hgroup}} \mult_{2 \eta_\Sb(y)}T_{[n]\setminus \Sb}f \ dy - \int_{\widehat{\Hgroup}} \mult_{2 \eta_\Sb(y)}T_{[n]\setminus \Sb}f \ dy = 0
	\end{align*}	
since $2 \eta_S(y)$ and $-2 \eta_S(y)$ are identically distributed. Therefore, $h\in L_p^\circ(\widehat{\Hgroup}; L_p(\mathcal{L}\Ggroup))$, so we raise to the power $p$ and average over $\Sb \subset [n]$ in \eqref{eq:lemma2Eq1}, and by \ref{itm:2} we get 
\begin{eqnarray*}
\lefteqn{\hskip-20pt \frac{1}{{n \choose k}} \suma{\substack{\Sb \subset [n] \\ |\Sb|=k}}  \int_{\widehat{\Hgroup}} \big\| \mult_{4m \eta_\Sb(y)} T_{[n]\setminus \Sb}f - T_{[n]\setminus \Sb}f \big\|_p^p \, dy} \\ \hskip20pt & \lesssim_p & \frac{m^p}{{n \choose k}} \suma{\substack{\Sb \subset [n] \\ |\Sb|=k}} \int_{\widehat{\Hgroup}} \big\| \cexp_{[n]\setminus \Sb} h(y) \big\|^p_p \, dy \\ & \lesssim_p & m^p \Bigg( \frac{k}{n} \Suma{j=1}{n} \int_{\widehat{\Hgroup}} \big\| (\partial_j \otimes \mathrm{id})h(y) \big\|_p^p \, dy + \Big(\frac{k}{n}\Big)^{\frac{p}{2}} \int_{\widehat{\Hgroup}} \|h(y)\|_p^p \, dy \Bigg) \\ [5pt] & \lesssim & m^p \Bigg( \frac{k}{n} \Suma{j=1}{n} \int_{\widehat{\Hgroup}} \big\| \partial_j \mult_{2 \eta(y)}f \big\|_p^p \, dy + \Big(\frac{k}{n}\Big)^{\frac{p}{2}} \int_{\widehat{\Hgroup}} \|\mult_{\eta(y)} f - f\|_p^p \, dy \Bigg).
\end{eqnarray*}
Therefore, putting altogether we get
\begin{eqnarray*}
\lefteqn{\frac{1}{{n \choose k}} \suma{\substack{\Sb \subset [n] \\ |\Sb|=k}} \int_{\widehat{\Hgroup}} \big\| \mult_{4m\eta_\Sb(y)}f - f \big\|_p^p \, dy} \\ & \lesssim & m^p \Bigg(\frac{k}{n} \Suma{j=1}{n} \int_{\widehat{\Hgroup}} \big\| \partial_j \mult_{2 \eta(y)}f \big\|^p_{p} \, dy + \Big(\Big(\frac{k}{n}\Big)^{\frac{p}{2}} + \frac{1}{m^p} \Big) \int_{\widehat{\Hgroup}} \big\| \mult_{\eta(y)}f - f \big\|^p_{p} \, dy \Bigg). 
\end{eqnarray*}
This completes the proof for any $m \geq \sqrt{n/k}$, as imposed in the statement. \fin

\begin{remark}
\emph{A close inspection of the proof of Theorem \ref{th:abstractThm} shows that one can exchange requirement \ref{itm:12} by the existence of a family of invertible operators $\{\mult_{y}\}_{y\in\widehat{\Hgroup}}$ satisfying the following properties:}
\begin{itemize}
\item Product structure. $\mult_y = \mult_{y_\Sb} \mult_{y_{[n]\setminus \Sb}}$.

\vskip8pt

\item Uniform boundedness. $\displaystyle \max_{|k|\leq 4m} \sup_{y \in \widehat{\Hgroup}} \big\| \mult_y^k \big\|_{p \to p} < \infty$.

\item Symmetry. \emph{$\mult_y^kf$ and $\mult_y^{-k}f$ have the same distribution function for each $f$.}
\end{itemize}
\emph{The above conditions are formally weaker than \ref{itm:12}, and more importantly, they clarify the fact that the auxiliary group $\HGroup$ is not necessary in the definition of $\Xp$-representable pairs. However, in the paper we have chosen to work with the stronger condition \ref{itm:12} because in all of our examples we are able to identify the associated group $\HGroup$ and this makes checking the conditions of Theorem \ref{th:abstractThm} simpler. }
\end{remark}

\begin{remark}
\emph{In many of the cases of interest there holds
\begin{equation}\label{eq:auxCondition}
\int_{\widehat{\Hgroup}} \big\| \partial_j \mult_{2\eta(y)} f \big\|_{L_p(\mathcal{L}(\Ggroup))}^p \, dy \lesssim_p \big\| \mult_{e_j}f - f \big\|_{L_p(\mathcal{L}(\Ggroup))}^p,
\end{equation}
for each $j=1, \ldots, n$. This improves the conclusion of Theorem \ref{th:abstractThm} to a truly metric inequality, which resembles the one in \cite{N16} more closely. This highlights the importance of the choice of the family of derivatives when checking conditions \ref{itm:2}-\ref{itm:12}. One case in which this holds is when derivatives and conditional expectations are related by
\begin{equation} \label{eq:derivativesExpectations}
	\partial_j = \mathrm{id} - \cexp_{\{j\}} \quad \mbox{for} \quad 1 \le j \le n.
\end{equation}
In turn, \eqref{eq:derivativesExpectations} holds when $\partial_j \bigchi_w(x) = \delta_{w_j \neq 0} \, \bigchi_w(x)$ for any character $\bigchi_{w}$ on $\widehat{\Hgroup}$.} 
\end{remark}


\section{\bf \large Intrinsic $\mathrm{X}_p$ inequalities} \label{sec2}

In this section we explore intrinsic $\mathrm{X}_p$ inequalities in the framework provided by Theorem \ref{th:abstractThm}. The term \lq\lq intrinsic" means here that there exists a natural inclusion of the index group $\Hgroup$ into the $\mathrm{X}_p$-group $\Ggroup$, so there is no need to use an auxiliary group $\HGroup$ other than $\Ggroup$. In our setting, this unfortunately forces us to work with $\Ggroup = \HGroup$ abelian. An extension of Theorem \ref{th:abstractThm} including noncommutative translations indexed by nonabelian groups $\Hgroup$ would open a door to potential applications in the metric geometry of noncommutative $L_p$-spaces. Intrinsic $\mathrm{X}_p$ inequalities will be complemented below with \lq\lq transferred $\mathrm{X}_p$ inequalities" which refer to those which admit some $\HGroup \neq \Ggroup$, including nonabelian $\mathrm{X}_p$-groups in the picture. We will illustrate this scenario in the context of free groups below. In the language of quantum group theory, intrinsic $\Hgroup$-translations are given by $\Ggroup$-comultiplications, while transferred $\Hgroup$-translations require more general $\Ggroup$-corepresentations in terms of $\HGroup$.  

\subsection{Continuous $\Xp$ inequalities} Let $\Hgroup=\Ggroup=\Z^n$, so that the dual groups are $n$-dimensional tori, which we identify with $[-1/2,1/2)^n$ for convenience. We also pick $\HGroup=\Z^n$, and we explore two different choices for $\eta:\mathbb{T}^n \to \mathbb{T}^n$. First, we take the map $\eta(y)_j = y_j/4m$, so each component of $\eta$ maps $\mathbb{T}^n$ to $[-1/(8m),1/(8m))$. The following encodes what follows from Theorem \ref{th:abstractThm} in this case. 

\begin{proposition}\label{prop:continuousXp}
If $p \ge 2$, $k \in [n]$ and $m \geq \sqrt{n/k}$, every $f : \mathbb{T}^n \rightarrow \C$ satisfies
\begin{eqnarray*}
\lefteqn{\frac{1}{{n \choose k}} \suma{\substack{\Sb \subset [n] \\ |\Sb|=k}} \int_{\mathbb{T}^n} \hskip-3pt \int_{\mathbb{T}^n} \big| f(x+ y_\Sb)-f(x) \big|^p \, dy  dx} \\ \!\!\!\! & \lesssim_p & \!\!\!\! m^p \hskip-3pt \int_{\mathbb{T}^n} \hskip-3pt \int_{\mathbb{T}^n} \hskip-3pt \Bigg\{ \frac{k}{n} \Suma{j=1}{n} \Big| f \big(x \hskip-1pt + \hskip-1pt \frac{y_j e_j}{4m}\big)-f(x)\Big|^p \hskip-3pt + \Big( \frac{k}{n} \Big)^{\frac{p}{2}} \Big| f\big(x \hskip-1pt + \hskip-1pt \frac{y}{4m}\big)-f(x)\Big|^p \Bigg\} dy dx. 
\end{eqnarray*}
\end{proposition}

\dem We start checking that $(\Z^n,\Z^n)$ is an $\Xp$-representable pair. Taking the group $\HGroup = \mathbb{Z}^n$,  it is clear that \ref{itm:12} holds with our choice of $\eta$ and multipliers given by the usual translations $\mathrm{M}_\gamma f(x) = f(x+\gamma)$ which are unitary modulations at the Fourier side. Second, if $\chi_w = \exp(2\pi i \langle w,\cdot\rangle)$ is the character associated with $w \in \Z^n$, we consider the differential operators
$$
\partial_j \chi_w = \delta_{w_j \neq 0} \, \chi_w \quad \mbox{for} \quad j \in [n].
$$
By \cite[Subsection 2.1]{CCP22}, we have
$$
\frac{1}{{n \choose k}} \suma{\substack{\Sb \subset [n] \\ |\Sb|=k}} \big\| \cexp_{[n]\setminus \Sb} f \big\|^p_{L_p(\mathbb{T}^n)} \lesssim_p \frac{k}{n} \Suma{j = 1}{n} \| \partial_j f \|_{L_p(\mathbb{T}^n)}^p + \Big(\frac{k}{n}\Big)^{\frac{p}{2}} \|f\|^p_{L_p(\mathbb{T}^n)}.
$$
In particular, by Fubini's theorem we get \ref{itm:2}. Then, Theorem \ref{th:abstractThm} yields
\begin{eqnarray} \label{eq:lemma9AltDerivatives}
\lefteqn{\hskip-5pt \frac{1}{{n \choose k}} \suma{\substack{\Sb \subset [n] \\ |\Sb|=k}} \int_{\mathbb{T}^n} \hskip-3pt \int_{\mathbb{T}^n}  \big| f(x+ y_\Sb) - f(x) \big|^p \, dy dx} \\ \nonumber \hskip30pt \!\!\!\! & \lesssim_p & \!\!\!\! m^p \int_{\mathbb{T}^n} \hskip-3pt \int_{\mathbb{T}^n} \frac{k}{n} \Suma{j=1}{n} \Big| \partial_{j}^y f\big(x+\frac{y}{2m}\big)\Big|^p + \Big( \frac{k}{n} \Big)^{\frac{p}{2}}  \Big| f\big(x+\frac{y}{4m}\big)-f(x)\Big|^p dy dx, 
\end{eqnarray}
using the superscripts in the partial derivatives to indicate the variable over which differentiation is performed. It only remains to estimate \eqref{eq:lemma9AltDerivatives} to establish the result. In order to do that, observe that \eqref{eq:derivativesExpectations} holds in this case. Applying it to $h(y) = f(x+y/2m)$ and using the properties of translations we get 
\begin{eqnarray*}
\lefteqn{\int_{\mathbb{T}^n} \hskip-3pt \int_{\mathbb{T}^n} \Big|\partial_{j}^y f\Big(x+\frac{y}{2m}\Big)\Big|^p dy dx} \\
	& = & \int_{\mathbb{T}^n} \hskip-3pt  \int_{\mathbb{T}^n} \Big|f\Big(x+\frac{y}{2m}\Big) - \cexp_{\{j\}}^y f\Big(x+\frac{y}{2m}\Big)\Big|^p dy dx \\ & = & \int_{\mathbb{T}^n} \hskip-3pt  \int_{\mathbb{T}^n} \Big|f\Big(x+\frac{y}{2m}\Big) - \int_{\mathbb{T}} f\Big(x+\frac{y_{[n]\setminus \{j\}}+te_j}{2m}\Big)dt \Big|^p dy dx \\ & \leq & \int_{\mathbb{T}^n} \hskip-3pt  \int_{\mathbb{T}^n} \hskip-3pt \int_{\mathbb{T}} \Big|f\Big(x+\frac{y}{2m}\Big) -  f\Big(x+\frac{y_{[n]\setminus \{j\}}+te_j}{2m}\Big) \Big|^p dt dy dx \\
	& = & \int_{\mathbb{T}^n} \hskip-3pt  \int_{\mathbb{T}^n} \hskip-3pt \int_{\mathbb{T}} \Big|f\Big(x+\frac{y_je_j}{2m}\Big) -  f\Big(x+\frac{te_j}{2m}\Big) \Big|^p dt dy dx \\ & = & \int_{\mathbb{T}^n} \hskip-3pt  \int_{\mathbb{T}^n} \Big|f\Big(x+\frac{y_je_j}{2m}\Big) -  f(x) \Big|^p dy dx \le 2 \int_{\mathbb{T}^n} \hskip-3pt  \int_{\mathbb{T}^n} \Big|f\Big(x+\frac{y_je_j}{4m}\Big) -  f(x) \Big|^p dy dx. 	
\end{eqnarray*}
Inserting that in the outcome of Theorem \ref{th:abstractThm} yields the assertion. \fin

\begin{remark}
\emph{A different choice of $\eta: \mathbb{T}^n\to\mathbb{T}^n$ is given by $\eta(y)_j = \sgn(y_j)/8m$ which leads to an inequality closer to \eqref{eq:metricXp}. Recalling that we identify $\mathbb{T}$ with $[-1/2,1/2)$, it turns out that $4m\eta(y) \equiv -e/2$ for $e=(1,1,\ldots,1)$, and one can prove the following statement for every $f:\mathbb{T}^n \rightarrow \C$
\begin{eqnarray} \label{eq:semicont}
\lefteqn{\hskip30pt \frac{1}{{n \choose k}} \suma{\substack{\Sb \subset [n] \\ |\Sb|=k}} \int_{\mathbb{T}^n} \Big| f\big(x+ \frac{e_\Sb}{2}\big)-f(x) \Big|^p  dx}  \\ \nonumber \!\!\!\! & \lesssim_p & \!\!\!\! m^p \int_{\mathbb{T}^n} \frac{k}{n} \Suma{j=1}{n} \Big|f\big(x+\frac{e_j}{8m}\big)-f(x)\Big|^p + \Big( \frac{k}{n} \Big)^{\frac{p}{2}} \frac{1}{2^n}\sum_{\varepsilon \in \Z_2^n} \Big|f\big(x+\frac{\varepsilon}{8m}\big)-f(x)\Big|^p dx. 
\end{eqnarray}
The proof of \eqref{eq:semicont} follows the same lines as that of Proposition \ref{prop:continuousXp}. It is also worth noting that inequality \eqref{eq:semicont} can be discretized to recover \eqref{eq:metricXp} and a limiting procedure allows one to pass from \eqref{eq:metricXp} to \eqref{eq:semicont}, we omit the details. This means that this \emph{semi-continuous} inequality is equivalent to Naor's original $\Xp$ inequality.} 
\end{remark}

\begin{remark}\label{th:metricXpTorus}
\emph{It is also possible to consider the usual differential structure on $\mathbb{T}^n$ which can be thought of as the most natural. In that case, $L_p$ valued versions of estimates for balanced truncations of Fourier series hold by Fubini's theorem \cite{CCP22} and one gets the following
\begin{eqnarray*}
\lefteqn{\frac{1}{{n \choose k}} \suma{\substack{\Sb \subset [n] \\ |\Sb|=k}} \int_{\mathbb{T}^n} \int_{\mathbb{T}^n} \big| f(x+y_\Sb)-f(x) \big|^p \, dy dx} \\ \!\!\!\! & \lesssim_p & \!\!\!\! m^p \int_{\mathbb{T}^n} \hskip-3pt \int_{\mathbb{T}^n} \frac{k}{n} \Suma{j=1}{n} \Big|\partial_{y_j} f\big(x+\frac{y}{2m}\big)\Big|^p + \Big( \frac{k}{n} \Big)^{\frac{p}{2}} \Big|f\big(x+\frac{y}{4m}\big)-f(x)\Big|^p dy dx,
\end{eqnarray*}
which can be seen to be weaker than those from Proposition \ref{prop:continuousXp}, see \cite{CCP22} for details. }
\end{remark}


\subsection{Cyclic groups with the word length} \label{subsec:2-2} We now focus on a discrete inequality that is a strict generalization of \eqref{eq:metricXp}. Consider the inclusion $\eta_{\ell} : \Z_{2\ell}^n \to \Z_{8\ell m}^n$ given by
\begin{align}\label{eq:mapEtaZ2l}
	\eta_{\ell}(x) = \big( \beta_\ell(x_1), \ldots,\beta_\ell(x_n) \big),
\end{align}
where $\beta_\ell(y) = y-\ell$ when $0 \leq y \leq \ell-1$ and $\beta_\ell(y)= y - (\ell-1)$ when $\ell \leq y \leq 2\ell-1$.

\begin{proposition}\label{th:metricXpCyclic}
If $p \ge 2$, $k \in [n]$ and $m \geq \sqrt{n/k}$, every $f : \Z_{8 \ell m}^n \rightarrow \C$ satisfies
\begin{eqnarray*}
\lefteqn{\hskip-60pt \frac{m^{-p}}{{n \choose k}} \suma{\substack{\Sb \subset [n] \\ |\Sb|=k}} \frac{1}{(8 \ell m)^n(2\ell)^n} \suma{x \in \Z_{8 \ell m}^n}  \suma{y \in \Z_{2\ell}^n} \big| f\big(x+4m\eta_{\ell}(y)_\Sb\big) - f(x)\big|^p} \\ \hskip40pt & \lesssim_p & (4\ell)^{p-1} \frac{k}{n} \Suma{j=1}{n} \frac{1}{(8\ell m)^n} \suma{x \in \Z_{8\ell m}^n} \big| f(x+e_j)-f(x) \big|^p \\ & + & \Big( \frac{k}{n} \Big)^{\frac{p}{2}} \frac{1}{(8 \ell m)^n(2\ell)^n} \suma{x \in \Z_{8 \ell m}^n} \suma{y \in \Z_{2\ell}^n} \big| f(x+\eta_{\ell}(y)) - f(x) \big|^p.
\end{eqnarray*} 
\end{proposition}

\dem Choose $\ell\geq 1$ and take $(\Hgroup,\Ggroup) = (\Z_{2\ell}^n,\Z_{8\ell m}^n)$ with auxiliary group $\HGroup = \Ggroup$. We want to check that this is an $\Xp$-representable pair. It is clear that the range of $\beta_\ell$ is $-[\ell] \cup [\ell] \subset \Z_{8\ell m}^n$. When $\ell=1$, it recovers the map that sends $0$ to $-1$ and $1$ to $1$. In addition, a simple observation yields that $-\beta_\ell(y) = \beta_\ell(2\ell-1- y)$ which implies condition \ref{itm:12} i). Next, denoting again characters by $\chi_w$ for each $w\in \Z_{2\ell}^n$, we choose conditional expectations given by
$$
\cexp_{[n] \setminus \Sb} \chi_w= \delta_{w \in \Z_{2\ell}^\Sb} \chi_w.
$$
By \cite[Subsection 2.2]{CCP22} and Fubini's theorem, \ref{itm:2} holds with derivatives given again by $\partial_j : \chi_w \mapsto \delta_{w_j \neq 0} \chi_w$. Finally, the map $\mult_\gamma : \chi_w \mapsto \chi_w(\gamma) \chi_w$ is a completely bounded multiplier on $L_p(\Z_{8\ell m}^n)$, so the family $\{\mult_\gamma\}_{\gamma \in \HGroup}$ satisfies \ref{itm:12} ii). The application of Theorem \ref{th:abstractThm} yields 
\begin{eqnarray} \label{eq:XpCyclicDeriv}
\lefteqn{\hskip-40pt \frac{m^{-p}}{{n \choose k}} \suma{\substack{\Sb \subset [n] \\ |\Sb|=k}} \frac{1}{(8 \ell m)^n(2\ell)^n} \suma{x \in \Z_{8 \ell m}^n}  \suma{y \in \Z_{2\ell}^n} \big| f\big(x+4m\eta_{\ell}(y)_\Sb\big) - f(x)\big|^p} \\ \nonumber \hskip40pt & \lesssim_p & \frac{k}{n} \Suma{j=1}{n} \frac{1}{(8\ell m)^n(2\ell)^n} \suma{x \in \Z_{8\ell m}^n} \suma{y \in \Z_{2\ell}^m} \big|\partial_{j}^yf \big(x+2\eta_{\ell}(y) \big) \big|^p \\ \nonumber & + & \Big( \frac{k}{n} \Big)^{\frac{p}{2}} \frac{1}{(8 \ell m)^n(2\ell)^n} \suma{x \in \Z_{8 \ell m}^n} \suma{y \in \Z_{2\ell}^n} \big| f(x+\eta_{\ell}(y)) - f(x) \big|^p.
\end{eqnarray} 
Since \eqref{eq:derivativesExpectations} holds in this setting too, we get 
\begin{align*}
		 \suma{x \in \Z_{8\ell m}^n} \suma{y \in \Z_{2\ell}^n}& \big| \partial_j^yf(x+2\eta_{\ell}(y)) \big|^p \\
		&= \suma{x \in \Z_{8\ell m}^n} \suma{y \in \Z_{2\ell}^n} \big| f(x+2\eta_{\ell}(y)) - \cexp_{\{j\}}^yf(x+2\eta_{\ell}(y)) \big|^p \\
		&= \suma{x \in \Z_{8\ell m}^n} \suma{y \in \Z_{2\ell}^n} \Big| \frac{1}{2\ell} \suma{t \in \Z_{2\ell}} f(x+2\eta_\ell(y)) - f(x+2\eta_\ell(y_{[n]\setminus\{j\}}+te_j)) \Big|^p \\
		&\leq \suma{x \in \Z_{8\ell m}^n} \suma{y \in \Z_{2\ell}^n} \frac{1}{2\ell} \suma{t \in \Z_{2\ell}} \big| f(x+2\eta_\ell(y)) - f(x+2\eta_\ell(y_{[n]\setminus\{j\}}+te_j)) \big|^p \\
		&= (2\ell)^{n-2} \suma{x \in \Z_{8\ell m}^n} \suma{y_j \in \Z_{2\ell}} \suma{t \in \Z_{2\ell}} \big| f(x+2 \beta_\ell(y_j)e_j) - f(x + 2 \beta_\ell(t)e_j) \big|^p,
\end{align*}
by the definition of $\cexp_{[n]\setminus \Sb}$ and convexity. Since $|\beta_\ell(y_j)-\beta_\ell(t)|\leq2\ell$, the translation invariance of the Haar measure yields
\begin{eqnarray*}
\lefteqn{\hskip-100pt (2\ell)^{n-2} \suma{x \in \Z_{8 \ell m}^n} \suma{y_j,t \in \Z_{2\ell}} \big| f(x + 2 \beta_\ell(y_j)e_j) - f(x+ 2 \beta_\ell(t)e_j) \big|^p} \\ \hskip70pt & \leq & (2\ell)^{n}(4\ell)^{p-1} \suma{x \in \Z_{8\ell m}^n} \big| f(x+e_j)-f(x) \big|^p,
\end{eqnarray*}
which yields the result when inserted in \eqref{eq:XpCyclicDeriv}. This completes the proof.  \fin

\begin{remark}
\emph{In  our cyclic metric $\Xp$ inequalities, a similar  computation as the one in \cite[Proposition 1.4]{NS16} shows that the condition $m \geq \sqrt{n/k}$ is necessary in Proposition \ref{th:metricXpCyclic}.}
\end{remark}

\section{\bf \large Transferred $\mathrm{X}_p$ inequalities}\label{sec4}

Now we apply Theorem \ref{th:abstractThm} to some pairs $(\Hgroup,\Ggroup)$ with $\Hgroup=\Z_{2}^n$ and nonabelian $\Ggroup$. Clearly, the choices of the auxiliary group $\HGroup$ and $\Ggroup$ in each case must be related. We give the details for one particular choice and comment on a few more later. We choose $\HGroup=\Z_{8m}^{n}$ and 
$$
\Ggroup= \mathbb{Z}_{8m}^{*n} := \Z_{8m} * \ldots * \Z_{8m}, 
$$
the free product of $n$ copies of $ \mathbb{Z}_{8m}$. For $j \in [n]$, denote the generator of its $j$-th factor by $g_j$. Elements $w \in \mathbb{Z}_{8m}^{*n}$ can and will be described as reduced words in the form
$$
w=g_{k_1}^{\ell_1} g_{k_2}^{\ell_2} \ldots g_{k_r}^{\ell_r} \quad \mbox{with} \quad k_i \neq k_{i+1} \quad \mbox{and} \quad \ell_i \in \mathbb{Z} \setminus 8m \mathbb{Z}.
$$
Given $u \in \Z_{8m}^{n}$, we consider the character $\bigchi_{u} : \Z_{8m}^{*n} \to \C$ determined by
$$\bigchi_{u}(w) = \exp \Bigg( \frac{2\pi i}{8m} \suma{j=1}^n u_j \suma{t  :  k_t = j} \ell_t \Bigg) \quad \mbox{for} \quad w = g_{k_1}^{\ell_1} g_{k_2}^{\ell_2} \ldots g_{k_r}^{\ell_r}.$$
It is clear that $\bigchi_{u}(w) \bigchi_{v}(w) = \bigchi_{u+v}(w)$. Next, we define the Fourier multipliers $\mult_u : \mathcal{L}(\Z_{8m}^{*n}) \rightarrow \mathcal{L}(\Z_{8m}^{*n})$ determined by $\lambda(w) \mapsto \bigchi_{u}(w)  \lambda(w)$ for the left regular representation $\lambda$ on $\Z_{8m}^{*n}$. 

\begin{lemma}\label{XpTransfMapMu}
$\mult_u$ is a normal unital trace-preserving $*-$homomorphism on $\mathcal{L}(\Z_{8m}^{*n})$.
\end{lemma}

\dem $\mult_u$ is clearly unital and
	\begin{equation*}
		\mult_{u}(\lambda(w)^*) = \mult_{u}(\lambda(w^{-1})) = \overline{\bigchi_{u}(w)} \ \lambda(w)^* = \mult_u (\lambda(w))^*.
	\end{equation*}
Next we claim that $\mult_u (\lambda(ww^\prime)) = \mult_u (\lambda(w)) \mult_u (\lambda(w^\prime))$. It is clear that this identity holds when both $w$ and $w'$ are powers of a fixed generator $g_j$. More precisely, we have  
	\begin{equation*}
		\mult_u (\lambda(g_j^{\ell_1+\ell_2})) = \mult_u (\lambda(g_j^{\ell_1})) \  \mult_u (\lambda(g_j^{\ell_2})).
	\end{equation*}
Now, let us consider two reduced words given by 
$$w = g_{k_1}^{\ell_1} g_{k_2}^{\ell_2} \ldots g_{k_r}^{\ell_r} \quad \mbox{and} \quad w^{\prime} = g_{k_1^\prime}^{\ell_1^\prime} g_{k_2^\prime}^{\ell_2^\prime} \ldots g_{k_s^\prime}^{\ell_s^\prime}.$$
If $k_r \neq k_1^\prime$ the claim trivially holds. If $k_r = k_1^\prime$ and $\ell_r + \ell_1^\prime \not\equiv 0 \ (\mbox{mod} \ 8m)$, then
		\begin{equation*}
			\mult_u(\lambda(ww^\prime)) = \mult_u(\lambda(g_{k_1}^{\ell_1}\ldots g_{k_{r-1}}^{\ell_{r-1}})) \ \mult_u(\lambda(g_{k_r}^{\ell_r+\ell_{1}^\prime})) \ \mult_u(\lambda(g_{k_2^\prime}^{\ell_2^\prime} \ldots g_{k_s^\prime}^{\ell_s^{\prime}}))
		\end{equation*}
which yields the same conclusion. For the remaining case, we may write $ww'$ as
\begin{equation*}
w w^\prime = \underbrace{g_{k_1}^{\ell_1} \ldots g_{k_{r-1}}^{\ell_{r-1}}}_{\rho} \underbrace{g_{k_2^\prime}^{\ell_2^\prime} \ldots g_{k_s^\prime}^{\ell_s^\prime}}_{\rho^\prime}.
\end{equation*}
Arguing as above, if $k_{r-1} \neq k'_2$ or $\ell_{r-1} + \ell_2' \not\equiv 0 \ (\mbox{mod} \ 8m)$ we get
\begin{eqnarray*}
\mult_u (\lambda(ww^\prime)) & = & \mult_u(\lambda(\rho)) \mult_u(\lambda(\rho^\prime)) \\ & = & \mult_u(\lambda(\rho)) \mult_u(\lambda(g_{k_r}^{\ell_r+\ell_1^\prime})) \mult_u(\lambda(\rho^\prime)) \ = \ \mult_u(\lambda(w)) \mult_u(\lambda(w^\prime)).
\end{eqnarray*}
One can iterate this process to deduce the claim. Thus, $\mult_u$ is a $*-$homomorphism on $\mathrm{span}\{\lambda(w):  w\in\Z_{8m}^{*n}\}$. In particular, it is a completely positive unital map. 
Moreover, it is trace-preserving since $\tau(\mult_u(\lambda(w))) = \bigchi_{u}(w) \tau(\lambda(w))$ and $\bigchi_{u}(e) = 1$. Finally, note that for any $f,g \in \mathcal{L}(\Z_{8m}^{*n})$ there holds
 \begin{eqnarray*}
 \tau(\mult_u(f)g^*) \!\!\!\! & = & \!\!\!\! \tau \Big[ \Big(\suma{w\in \mathbb{Z}_{8m}^{*n}} \widehat{f}(w\ \bigchi_{u}(w\ \lambda(w) \Big) \Big( \suma{\eta \in \mathbb{Z}_{8m}^{*n}} \widehat{g}(\eta) \ \lambda(\eta)\Big)^* \Big] \\ \!\!\!\! & = & \!\!\!\! \suma{w \in \mathbb{Z}_{8m}^{*n}} \widehat{f}(w) \bigchi_{u}(w) \overline{\widehat{g}(w)} \, = \suma{w \in \mathbb{Z}_{8m}^{*n}} \widehat{f}(w) \overline{\bigchi_{-u}(w)} \overline{\widehat{g}(w)} \, = \, \tau(f (\mult_{-u}g)^*).
 \end{eqnarray*}
 Since $\mult_{-u}$ extends to a bounded map on $L_1(\mathcal{L}(\Z_{8m}^{*n}))$, $\mult_u$ is $w^*-$continuous. \fin 
 
Before exploring $\mathrm{X}_p$ inequalities for free groups, we shall need to use Theorem B in \cite{CCP22} with values in a noncommutative $L_p$-space over a QWEP von Neumann algebra $\mathcal{M}$. More precisely, given a mean-zero $f: \Z_2^n \to \mathcal{M}$ we claim that 
\begin{eqnarray} \label{opvalued}
\lefteqn{\hskip-30pt \frac{1}{{n \choose k}} \suma{\substack{\Sb \subset [n] \\ |\Sb|=k}} \big\|(\cexp_{[n]\setminus \Sb} \otimes \mathrm{id}) f \big\|^p_{L_p(\Z_2^n;L_p(\mathcal{M}))}} \\ \nonumber \hskip40pt & \lesssim_p & \frac{k}{n} \Suma{j=1}{n} \big\|(\partial_j \otimes \mathrm{id}) f \big\|^p_{L_p(\Z_2^n;L_p(\mathcal{M}))} + \Big(\frac{k}{n}\Big)^{\frac{p}{2}} \|f\|^p_{L_p(\Z_2^n;L_p(\mathcal{M}))}, 
\end{eqnarray}
where $\cexp_{[n]\setminus \Sb}$ and $\partial_j$ stand for the original conditional expectation and directional derivative used by Naor in \cite{N16}. To check that \eqref{opvalued} holds, one can see that the proof in \cite{CCP22} goes through in the operator valued setting as long as dimension-free bounds hold for operator valued Riesz transforms \cite{JMP18}, which just requires to use the same argument together with Fubini's theorem from \cite{J04}. 

\begin{proposition} \label{metricXpTransfer}
If $p \ge 2$, $k \in [n]$ and $m \geq \sqrt{n/k}$, every $f \in L_p(\mathcal{L}(\Z_{8m}^{*n}))$ satisfies
\begin{eqnarray*}
\lefteqn{\hskip-40pt \frac{1}{{n \choose k}} \suma{\substack{\Sb \subset [n] \\ |\Sb|=k}} \frac{1}{2^n} \suma{\varepsilon \in \Z_2^n} \big\|\mult_{4m\varepsilon_\Sb}f - f \big\|^p_{L_p(\mathcal{L}(\Z_{8m}^{*n}))}} \\ \hskip40pt & \lesssim_p & m^p \Bigg( \frac{k}{n} \Suma{j=1}{n} \big\| \mult_{e_j}f - f \big\|_p^p + \Big(\frac{k}{n} \Big)^{\frac{p}{2}} \frac{1}{2^n} \suma{\varepsilon \in \Z_2^n} \big\| \mult_{\varepsilon}f - f \big\|_p^p \Bigg).
\end{eqnarray*}
\end{proposition}

\dem We apply Theorem \ref{th:abstractThm} to $(\Hgroup, \Ggroup) = (\Z_2^n,\Z_{8m}^{*n})$, so we need to check that it is an $\Xp$-representable pair. Condition \ref{itm:2} follows from \eqref{opvalued} since $\mathcal{L}(\Z_{8m}^{*n})$ is a QWEP von Neumann algebra. Consider the auxiliary group $\HGroup = \mathbb{Z}_{8m}^n$ with the inclusion $\eta$ defined as in \eqref{eq:mapEtaZ2l} for $\ell=1$. Condition \ref{itm:12} i) is then automatically satisfied. Finally, according to Lemma \ref{XpTransfMapMu} we know that $\mult_\gamma$ extends to a completely contractive map on $L_p(\mathcal{L}(\Z_{8m}^{*n}))$ for any $1 \leq p \leq \infty$. Therefore, the family $\{\mult_\gamma\}_{\gamma \in\Z_{8m}^n}$ can be used to check \ref{itm:12} ii) and Theorem \ref{th:abstractThm} leads to
\begin{align*}
\frac{1}{{n \choose k}} \suma{\substack{\Sb \subset [n] \\ |\Sb|=k}}& \frac{1}{2^n} \suma{\varepsilon\in\Z_2^n} \big\|\mult_{4m\varepsilon_\Sb}f - f \big\|^p_{L_p(\mathcal{L}(\Z_{8m}^{*n}))}\\
	&\lesssim_p m^p\Bigg(\frac{k}{n} \Suma{j=1}{n} \frac{1}{2^n} \sum_{\varepsilon\in\Z_2^n} \big\|(\partial_j \otimes \mathrm{id})\mult_{2\varepsilon}f \big\|^p_p + \Big(\frac{k}{n} \Big)^{\frac{p}{2}} \frac{1}{2^n} \suma{\varepsilon\in\Z_2^n} \big\|\mult_{\varepsilon}f - f \big\|_p^p\Bigg),
\end{align*}
for every $f \in L_p(\Z_{8m}^n; L_p(\Z_{8m}^{*n}))$. Finally, the computations for discrete derivatives from Subsection~\ref{subsec:2-2} imply that condition \eqref{eq:auxCondition} holds, yielding the result. \fin

The example above is clearly not the only possible choice that we can make for the groups $\HGroup$ and $\Ggroup$ in the pair. In general, given abelian groups $\HGroup_1, \ldots , \HGroup_n$ and $\HGroup = \HGroup_1 \times \HGroup \times \cdots \times \HGroup_n$ so that there is a nice map
$$
\eta: \Z_2^n \to \widehat{\HGroup},
$$
one can take $\Ggroup = \HGroup_1 * \ldots * \HGroup_n$ and construct a family of multipliers $\{\mult_\gamma\}_{\gamma \in \widehat{\HGroup}}$ as above. Once the properties of the family are checked, the above scheme yields a metric $\Xp$ inequality for the pair $(\Z_2^n,\HGroup,\Ggroup)$. For example, the interested reader can check that a representation in the spirit of Lemma \ref{XpTransfMapMu} for $\HGroup=\Z^n$ and $\Ggroup = \mathbb{F}_n$ (the free group of $n$ generators) holds. This yields the following free metric $\Xp$ inequality in the free group algebra
\begin{align*}
	\frac{1}{{n \choose k}}&  \suma{\substack{\Sb \subset [n] \\ |\Sb|=k}} \frac{1}{2^n} \suma{\varepsilon \in \mathbb{Z}_{2}^n} \big\|\mult_{4m\varepsilon_\Sb}f - f\big\|^p_{L_p(\mathcal{L}(\mathbb{F}_n))}\\
	&\lesssim_p m^p \Bigg( \frac{k}{n} \Suma{j=1}{n} \big\|\mult_{e_j}f -f 
	\big\|^p_{L_p(\mathcal{L}(\mathbb{F}_n))} + \Big(\frac{k}{n} \Big)^{\frac{p}{2}} \frac{1}{2^n} \suma{\varepsilon \in \mathbb{Z}_{2}^n} \big\|\mult_{\varepsilon}f - f \big\|_{L_p(\mathcal{L}(\mathbb{F}_n))}^p \Bigg).
\end{align*}


\section{\bf \large Metric consequences} 

We conclude by collecting a few metric consequences of the results in the previous sections. In \cite{N16}, Naor gave two sufficient hypothesis for a Banach space $\mathbb{X}$ to be a metric $\Xp$ space: a Banach $\Xp$ inequality and operator valued dimension-free estimates for Riesz transforms. In particular, any mean-zero $f:\mathbb{Z}_2^n\to \mathbb{X}$ must satisfy
$$
\frac{1}{2^n} \sum_{\varepsilon \in \mathbb{Z}_2^n} \Big\| \Suma{j=1}{n} \varepsilon_j (\partial_j \otimes \mathrm{id})f \Big\|^p_{L_p(\mathbb{Z}_2^n;\mathbb{X})} \simeq_p \big\|(\Delta^{-1/2} \otimes \mathrm{id})f\big\|^p_{L_p(\mathbb{Z}_2^n;\mathbb{X})}.
$$
When $\mathbb{X}=L_p(\mathcal{M})$, the Banach $\Xp$ inequality is \cite[Theorem 1.2]{CCP22}, while the Riesz transform estimates can be found in \cite{JMP18} when $\mathcal{M}$ is QWEP, as explained in the previous section. This latter fact had previously been announced for the Schatten classes $S_p$ in \cite{N16,NPS20}. The result that we get is the following:

\begin{theorem}[Noncommutative $L_p$ spaces are metric $\Xp$ spaces]\label{thm:metricXpncLp}
Let $\mathcal{M}$ be a \emph{QWEP} von Neuman algebra. Then, if $p \ge 2$, $k \in [n]$ and $m \geq \sqrt{n/k}$, every $f : \Z_{8m}^n \rightarrow L_p(\mathcal{M})$ satisfies
\begin{eqnarray*}
\lefteqn{\hskip-60pt \frac{m^{-p}}{{n \choose k}} \suma{\substack{\Sb \subset [n] \\ |\Sb|=k}} \frac{1}{(16m)^n} \suma{x \in \Z_{8m}^n} \suma{\varepsilon \in \mathbb{Z}_2^n} \big\|f(x+4m\varepsilon_\Sb)-f(x)\big\|^p_{L_p(\mathcal{M})}} \\ \hskip60pt & \lesssim_p & \frac{k}{n} \Suma{j=1}{n} \frac{1}{(8m)^n} \suma{x \in \Z_{8m}^n} \big\|f(x+e_j)-f(x)\big\|^p_{L_p(\mathcal{M})} \\ & + & \Big(\frac{k}{n}\Big)^{\frac{p}{2}} \frac{1}{(16m)^n} \suma{x\in\Z_{8m}^n} \suma{\varepsilon \in \mathbb{Z}_2^n} \big\|f(x+\varepsilon)-f(x)\big\|^p_{L_p(\mathcal{M})}.
\end{eqnarray*}
\end{theorem}

Theorem \ref{thm:metricXpncLp} provides a large class of examples of metric $\Xp$ spaces beyond the class of (commutative) $L_p$ spaces. Also, as shown in \cite{N16}, being a metric $\Xp$ space implies lower estimates on the distorsion of bi-Lipschitz embeddings of nonlinear sets. Indeed, let $c_{L_p(\mathcal{M})}(\mathbf{X})$ denote the smallest norm of a bi-Lipschitz map $\mathbf{X} \to L_p(\mathcal{M})$. Then we have the following:  

\begin{corollary} \label{CorollaryDistortion}
Let $\mathcal{M}$ be \emph{QWEP} and $2 < q < p$. Then
\begin{itemize}
\item[\emph{i)}] $\displaystyle c_{L_p(\mathcal{M})} \big( [m]_q^n \big) \, \asymp_{p,q} \, \min \Big\{n^{\frac{(p-q)(q-2)}{q^2(p-2)}}, m^{1 - \frac{2}{q}} \Big\}$.

\vskip3pt

\item[\emph{ii)}] If $c_{L_p(\mathcal{M})} \big((L_q,\|x-y\|_q^\theta) \big) < \infty$, then necessarily $\theta \leq q/p$.
\end{itemize}
\end{corollary}

We refer to \cite{N16,NS16} for precise definitions on $[m]_q^n$ and $\theta$-snowflakes. We end with a note about the optimality of the discrete cyclic inequalities from Section \ref{sec2}. Given $m,n \geq 2$, we know from \cite[Lemma 3.1]{NS16} that there exists $h_m^n : \Z_{m}^n \rightarrow \{0,\ldots,4m\}^{2n}$ such that 
$$\Bigg( \Suma{j=1}{n} \Big| \exp \Big( \frac{2\pi i x_j}{m} \Big) - \exp \Big(\frac{2\pi i y_j}{m} \Big) \Big|^q \Bigg)^{\frac{1}{q}} \, \sim \, \big\| h_m^n(x) - h_m^n(y) \big\|_q$$
holds (up to absolute constants) for any $x,y \in \Z_{m}^n$ and $q \ge 2$. Given a bi-Lipschitz map $g : [16m]^{2n}_q \rightarrow L_p$ with bi-Lipschitz norm $D$, set $F= g \circ h_{4m,n}$. Then, arguing as in \cite[Theorem 1.14]{NS16} we get  
$$\frac{1}{{n \choose k}} \suma{\substack{\Sb \subset [n] \\ |\Sb|=k}} \frac{1}{(4m)^{2n}} \suma{x,y \in \Z_{4m}^{n}} \big\|F(x+y_\Sb)-F(x) \big\|_{L_p}^p \, \gtrsim_p \, m^p k^{\frac{p}{q}} \big|e^{\frac{2\pi i}{4m}}-1\big|.$$
Next, we make the choices $k = \lceil n^{\frac{p(q-2)}{q(p-2)}} \rceil $ and $m = \lceil n^{\frac{p-q}{q(p-2)}} \rceil$, which ensure that $m \gtrsim \sqrt{n/k}$, and so we can apply Proposition \ref{th:metricXpCyclic} to $F$ in order to get the following statement: if $2<q<p$ and $m,n, \ell \in \N$, there holds
$$
c_{L_p}([m]_q^n) \gtrsim_{p,q} \min\Big\{ |e^{i \pi/\ell}-1| \, n^{\frac{(p-q)(q-2)}{q^2(p-2)}}, m^{1 - \frac{2}{q}} \Big\} \gtrsim_{p,q,\ell} \min\Big\{ n^{\frac{(p-q)(q-2)}{q^2(p-2)}}, m^{1 - \frac{2}{q}} \Big\}.
$$
Thus, up to constants (which get worse as $\ell \to \infty$) we recover the optimal distortions found in \cite{N16}. Also,  Naor's result for $\theta$-snowflakes in \cite{N16} follow from Proposition \ref{th:metricXpCyclic}. 

\begin{problem}
According to \cite[Theorem 12]{NPS20} 
\begin{equation} \label{EqNPS}
c_{S_p}(\mathrm{X}) \le \dim(\mathrm{X})^{\frac{q}{2}(\frac1q - \frac1p)}
\end{equation}
for every linear subspace $\mathrm{X} \subset S_q^n$ and any $2 < q < p < \infty$. Is there an embedding of $\mathrm{X}$ into $S_p(\ell_q^m)$ with $m = \dim(\mathrm{X})$ and constants independent of $m$? This is trivially true for row, column or diagonal subspaces and we know from \cite[Theorem 2]{JP} that it holds for the full space $\mathrm{X}=S_q^n$ with optimal $m$ being $\dim(\mathrm{X}) = n^2$. If this was the case for general $\mathrm{X}$, we would get $c_{S_p}(\mathrm{X}) \lesssim c_{S_p}(S_p(\ell_q^m))$ whose rightmost term is dominated by an operator space $($matrix amplification$)$ form of $c_{S_p}(\ell_q^m)$. Using that $S_p$ is a metric $\mathrm{X}_p$-space and arguing as above 
\begin{equation} \label{EqNaor}
c_{S_p}(\ell_q^m) \asymp_{p,q} m^{\frac{(p-q)(q-2)}{q^2(p-2)}}.
\end{equation}
Does $c_{S_p}(S_p(\ell_q^m))$ behave like $c_{S_p}(\ell_q^m)$? If so, one could guess that 
\begin{equation} \label{EqConjecture}
c_{S_p}(\mathrm{X}) \le_{p,q} \dim(\mathrm{X})^{\frac{(p-q)(q-2)}{q^2(p-2)}} \quad \mbox{for every} \quad \mathrm{X} \subset S_q^n?
\end{equation}
This would already improve \eqref{EqNPS} and it is best possible for $\mathrm{X} = \ell_q^n = \mathrm{diag}(S_q^n)$. It is still open to find the optimal distortion $c_{S_p}(\mathrm{X})$ for the full space $\mathrm{X} = S_q^n$. In this case, we have the following bounds 
\begin{equation} \label{EqTrivial}
\mathrm{A} := n^{\frac{(p-q)(q-2)}{q^2(p-2)}} \le c_{S_p}(S_q^n) \le \min \Big\{ n^{\frac12 - \frac1q}, n^{\frac1q - \frac1p} \Big\}.
\end{equation}   
Indeed, the lower bound follows from \eqref{EqNaor} and the inclusion $\ell_q^n \subset S_q^n$. The upper bound follows easily using H\"older inequality. Note that \eqref{EqConjecture} gives $c_{S_p}(S_q^n) \le \mathrm{A}^2$ but \eqref{EqTrivial} implies that no bound $c_{S_p}(S_q^n) \le \mathrm{A}^\beta$ could be optimal for any $\beta > 1$.     
\end{problem}

\subsubsection*{\bf Acknowledgement.} 

The third-named author wants to express his gratitude to Alexandros Eskenazis and Assaf Naor for \hskip-1pt some \hskip-1pt valuable comments on a preliminary version, which led to an improved presentation. The three authors were partly supported by the Spanish Grant PID 2019-107914GB-I00 (MCIN / PI J. Parcet) as well as Severo Ochoa Grant CEX2019-000904-S (ICMAT), funded by MCIN/AEI 10.13039/501100011033. Jos\'e M. Conde-Alonso was also supported by the Madrid Government Program V under PRICIT Grant SI1/PJI/2019-00514. Antonio Ismael Cano-Mármol was supported by the Grant SEV-2015-0554-19-3 which is funded by MCIN/AEI/10.13039/501100011033 and by ESF Investing in your future.

\bibliography{biblio}
\bibliographystyle{plain}

\hfill \noindent \textbf{Antonio Ismael Cano M\'armol} \\
\null \hfill Instituto de Ciencias Matem{\'a}ticas 
\\ \null \hfill Consejo Superior de Investigaciones Cient{\'\i}ficas 
\\ \null \hfill\texttt{ismael.cano@icmat.es}

\vskip2pt

\hfill \noindent \textbf{Jos\'e M. Conde-Alonso} \\
\null \hfill Universidad Aut\'onoma de Madrid 
\\ \null \hfill Instituto de Ciencias Matem{\'a}ticas  
\\ \null \hfill\texttt{jose.conde@uam.es}

\vskip2pt

\hfill \noindent \textbf{Javier Parcet} \\
\null \hfill Instituto de Ciencias Matem{\'a}ticas 
\\ \null \hfill Consejo Superior de Investigaciones Cient{\'\i}ficas 
\\ \null \hfill\texttt{parcet@icmat.es}

\end{document}